\documentclass[journal]{IEEEtran}
\usepackage{ifpdf}
\usepackage{cite}
\ifCLASSINFOpdf
\usepackage[pdftex]{graphicx}
\usepackage{epstopdf}
\newcommand*{\Scale}[2][4]{\scalebox{#1}{\ensuremath{#2}}}%
\else
\usepackage[dvips]{graphicx}
\newcommand*{\Scale}[2][4]{\scalebox{#1}{\ensuremath{#2}}}%
\fi
\usepackage[cmex10]{amsmath}
\usepackage{mathtools}                  %
\usepackage{amssymb}					%Calls AMS symbols
\usepackage{arydshln}
\usepackage{amsfonts}
\usepackage{blkarray}

\usepackage{algorithmic}
\usepackage{array}
\usepackage{cases}
\usepackage{mdwmath}
\usepackage{mdwtab}
\usepackage{amsthm}
\usepackage{tikz}

\newcommand*\circled[1]{\tikz[baseline=(char.base)]{ \node[shape=circle,draw,inner sep=1pt] (char) {#1};}}

\newcommand\hlight[1]{\tikz[overlay, remember picture,baseline=-\the\dimexpr\fontdimen22\textfont2\relax]\node[rectangle,fill=blue!50,rounded corners,fill opacity = 0.2,draw,thick,text opacity =1] {$#1$};} 

\makeatletter
\renewcommand*\env@matrix[1][*\c@MaxMatrixCols c]{%
	\hskip -\arraycolsep
	\let\@ifnextchar\new@ifnextchar
	\array{#1}}
\makeatother

\usepackage[tight,footnotesize]{subfigure}

\usepackage{stfloats}
\usepackage{multirow}

\begin{document}
% at, but, by, for, in, nor, of, on, or, the, to and up, which are usually	% not capitalized unless they are the first or last word of the title.
\title{Sparse Tableau Formulation for Optimal Power Flow Applications}
	
\author{Byungkwon~Park,~\IEEEmembership{Student Member,~IEEE,}
		~Jayanth Netha,
		~Michael C. Ferris,
		and~Christopher L. DeMarco,~\IEEEmembership{Member,~IEEE}% <-this % stops a space
\thanks{B. Park and C. L. DeMarco are with the Department of Electrical and Computer Engineering, University of Wisconsin-Madison, Madison, WI 53706 USA e-mail:(bpark52@wisc.edu; cdemarco@wisc.edu.}% <-this % stops a space
\thanks{J. Netha is with the Department of Industrial Systems and Engineering, University of Wisconsin-Madison, Madison, WI 53706 USA (e-mail: jnetha@wisc.edu).}% <-this % stops a space
\thanks{M. C. Ferris is with the Department of Computer Sciences, University of
	Wisconsin-Madison, Madison, WI 53706 USA (e-mail: ferris@cs.wisc.edu).}
}
	

\markboth{}%
{Shell \MakeLowercase{\textit{et al.}}: Bare Demo of IEEEtran.cls for IEEE Journals}



\maketitle
	
\begin{abstract}
Typical formulations of the optimal power flow (OPF) problem rely on what is termed the ``bus-branch'' model, with network electrical behavior summarized in the $Y_{bus}$ admittance matrix. From a circuit perspective, this admittance representation restricts network elements to be voltage controlled and  limitations of the $Y_{bus}$ have long been recognized. A fixed $Y_{bus}$ is unable to represent an ideal circuit breaker, and more subtle limitations appear in transformer modeling. In power systems parlance, more detailed approaches to overcome these limitations are termed ``node-breaker'' representations, but these are often cumbersome, and are not widely utilized in OPF. This paper develops a general network representation adapted to the needs of OPF, based on the Sparse Tableau Formulation (STF) with following advantages for OPF: (i) conceptual clarity in formulating constraints, allowing a comprehensive set of network electrical variables; (ii) improved fidelity in capturing physical behavior and engineering limits; (iii) added flexibility in optimization solution, in that elimination of intermediate variables is left to the optimization algorithm. The STF is then applied to OPF numerical case studies which demonstrate that the STF shows little or no penalty in computational speed compared to classic OPF representations, and sometimes provides considerable advantage in computational speed.
\end{abstract}
	

% \fi
\IEEEpeerreviewmaketitle

\section{Introduction}
	
\IEEEPARstart{U}{}nder the banner of the Smart Grid \cite{KassakianSchalensee2011}, power systems today see growing integration of technologies from communications, advanced control, signal processing, power electronics, and data analytics, coupled with improving efficiency and cost effectiveness of  distributed energy resources. These trends open the door to a much wider range of control actions and decision variables in grid planning and operation, and motivate new optimization approaches to exploit these opportunities.

One of the central optimization problem underlying grid planning and operation is optimal power flow (OPF) \cite{MaryB.Cain2012}. The vast majority of OPF  approaches  formulate the  network constraints  based on the bus admittance matrix $Y_{bus}$, in which the bus voltage phasors serve as the key ``state'' variables, analogous to a ``strict" nodal analysis in standard circuit theory \cite{SadikuAlexander2011}, \cite{Grainger1994}. However, a strict nodal analysis disallows many standard circuit elements by its requirement that each element's current(s) be expressible a  function of its current(s) \cite{ChuaDesoerKuh1987}. In a power systems context, the $Y_{bus}$ formulation imposes similar restrictions; e.g. a fixed $Y_{bus}$ is unable to represent ideal circuit breakers in the network, because one cannot describe the current through the element as a function of voltage when the breaker is closed. These limitations spur growing recognition of the value of node-breaker representations \cite{GeneralElectric2016}, \cite{ThomasKincicDaviesEtAl2016}, that allow realistic representation of substation reconfiguration via circuit breakers in contingency analysis. However, much of the literature seeking to develop advanced OPF algorithms has remained focused on $Y_{bus}$ formulations, and lacks the generality of node-breaker representations.

The work of this paper will seek to gain the flexibility of general node-breaker formulations, while adopting a straightforward, algorithmic approach to network constraint formulation that is well suited to the OPF. To this end, it is worthwhile to briefly review comparable developments in the history of computer-aided analysis tools for electronic circuit design. While this is a vast literature, relevant early milestones applying optimization in automated network design include \cite{FischlPuntel1972}, \cite{W.R.Punteletal.1973} and in particular \cite{HACHTELBRAYTONGUSTAVSON1972}. The ``\textbf{Sparse Tableau Formulation} (STF)'' was particularly advocated by IBM for electronic design in the context of circuit optimization. While it failed to achieve the wide-spread adoption enjoyed by its contemporary and competitor SPICE \cite{NagelPederson1973}, the circuit analysis program ASTAP \cite{IBMPPDSH1973},\cite{WEEKSJIMENEZMAHONEYEtAl1973} developed by IBM successfully utilized STF. The  benefits of a STF-based formulation for power flow equations were explored in the late 1970's by \cite{DlRECTORSULLIVAN1979}, but have received little attention in subsequent decades. With advances in optimization algorithms, and in particular with automated elimination techniques that reduce penalties associated with the retention of large numbers of variables in a sparse formulation, this paper seeks to demonstrate that the STF is particularly well suited to Optimal Power Flow.

The organization of this paper is as follow. Section \ref{chap:3.2} reviews the general Sparse Tableau Formulation from a standard circuit analysis perspective. Section \ref{chap:3.3} then examines special cases and requirements associated with the power system application that allow simplification of general Sparse Tableau Formulation. Section \ref{chap:3.4} discusses the relationship between the STF and the $Y_{bus}$ in those cases for which both may be applied, examines non-typical network elements for whose modelling the Sparse Tableau Formulation provides particuar advantageous. The application to representative OPF examples in a general purpose optimization tool \cite{GAMS}, along with comparisons of computational speed between the STF and $Y_{bus}$ formulations, is described in Section \ref{chap:3.5}.
	
\section{Background}
\label{chap:3.2}
Here, we review the key steps in constructing STF circuit constraint equations \cite{ChuaDesoerKuh1987}.  As case of  interest in power systems, we give special attention to two-port circuit elements, and assume that the circuit analysis is conducted with respect to complex phasor branch or port currents, denoted $i$, complex phasor branch or port voltages, denoted $v$, and complex node voltages, denoted $V$. Note  that all equality constraints below are complex.

\noindent\textbf{Step 1.} Write a complete set of linearly independent KCL equations, employing the node-to-element reduced incidence matrix $A$: 
\begin{equation}
	\label{eq:KCL}	
	Ai = 0
\end{equation}
\textbf{Step 2.} Write a complete set of linearly independent KVL equations:
\begin{equation}
	\label{eq:KVL}	
	v-A^TV = 0
\end{equation}
\textbf{Step 3.} Write the element constituative equations. If the circuit elements are all affine linear, these may be written as:
\begin{equation}
\label{eq:branch}	
	F_vv + F_ii = u_s
\end{equation}

Equations (\ref{eq:KCL}), (\ref{eq:KVL}), and (\ref{eq:branch}) are the tableau equations. For an element represented as two-port, $v$ and $i$ quantities appear in ``port-pairs," as illustrated in the figure bleow. 
\begin{figure}[!htb]
	\centering
	\includegraphics[height=1in,width=2in]{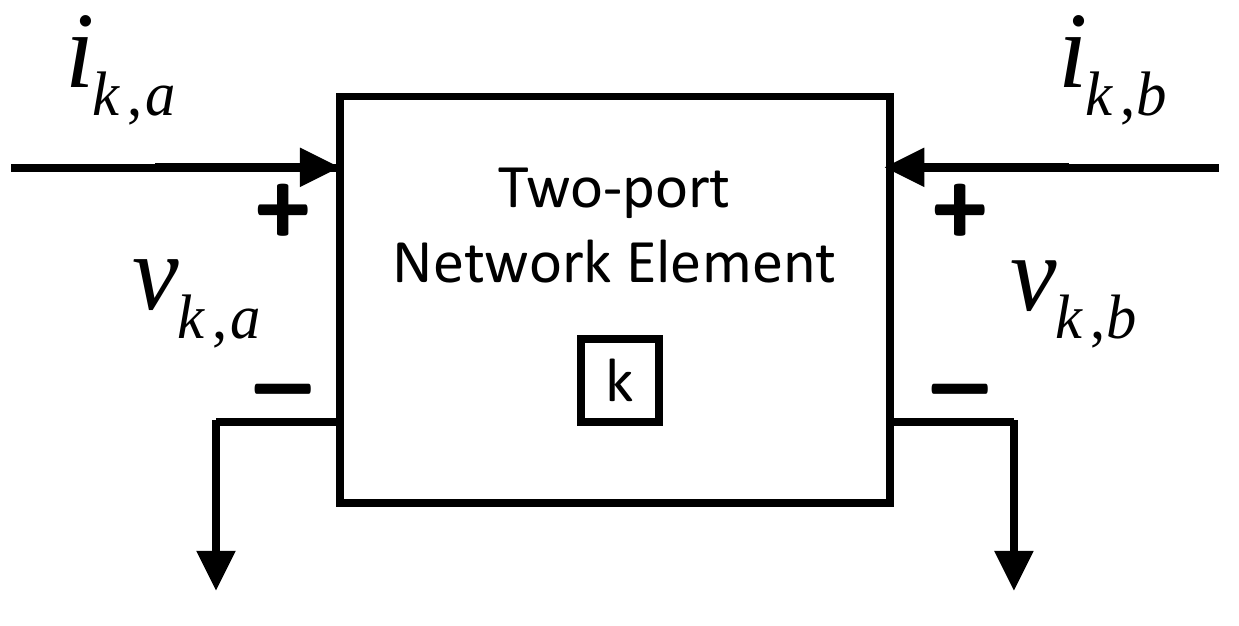}
	\caption{Two port representation and its matrix form}
	\label{fig:twoport}
\end{figure}

A transmisson line is the two-port element appearing perhaps most commonly in the power systems model. In this context, port $a$ quantites are typically  termed ``sending end" positive sequence voltage and current, and port $b$  quantities termed  ``receiving end."  A two-port element's constitutive relations place two independent algebraic constraints on the four variables $(v_a, i_a, v_b, i_b)$ to specify the element's behavior.  For each network element $k$, in the case of complex, phasor-based analysis,  these constraints take the general implicit form   
\begin{align}
\label{eq:general_element}
f_k : \mathbb{C}^{4} \longrightarrow \mathbb{C}^{2} \\ \nonumber
f_k (v_{k,a}, i_{k,a}, v_{k,b}, i_{k,b}) =
\begin{bmatrix}
0 \\
0 \\ 
\end{bmatrix} 
\end{align}

With all element equations composed together as $f(v, i) = \mathbf{0}$, observe that the element constitutive equation (\ref{eq:branch}) in the linear tableau equation is a special affine case of $f(v, i) = \mathbf{0}$;  in particular:
\begin{align}
f(v, i) = \mathbf{0}  \triangleq 
F_{v}v + F_{i}i - u_s = 0
\end{align}

\begin{figure}[htb!]
	\centering
	\subfigure[Linear element circuit]
	{\includegraphics[height=1.1in, width=2.5in]{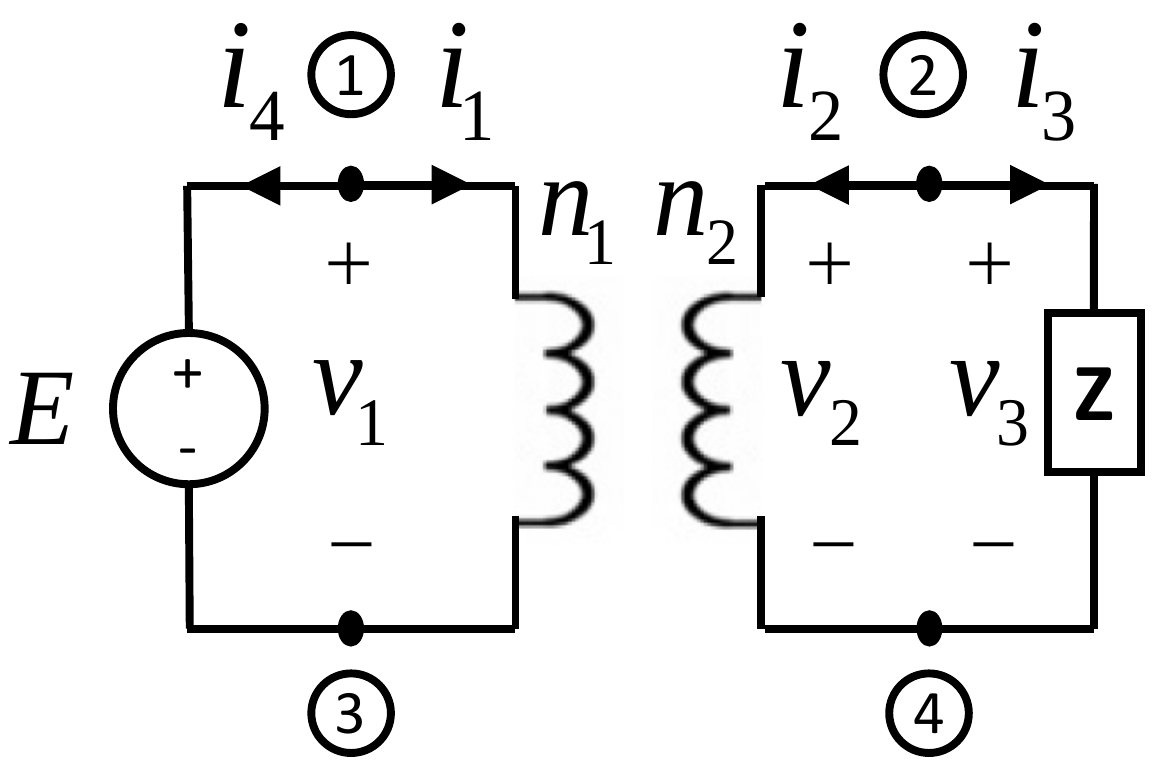}
		\label{fig:first_case}}
	\hfil
	\subfigure[Connected digraph]
	{\includegraphics[height=0.9in, width=2in]{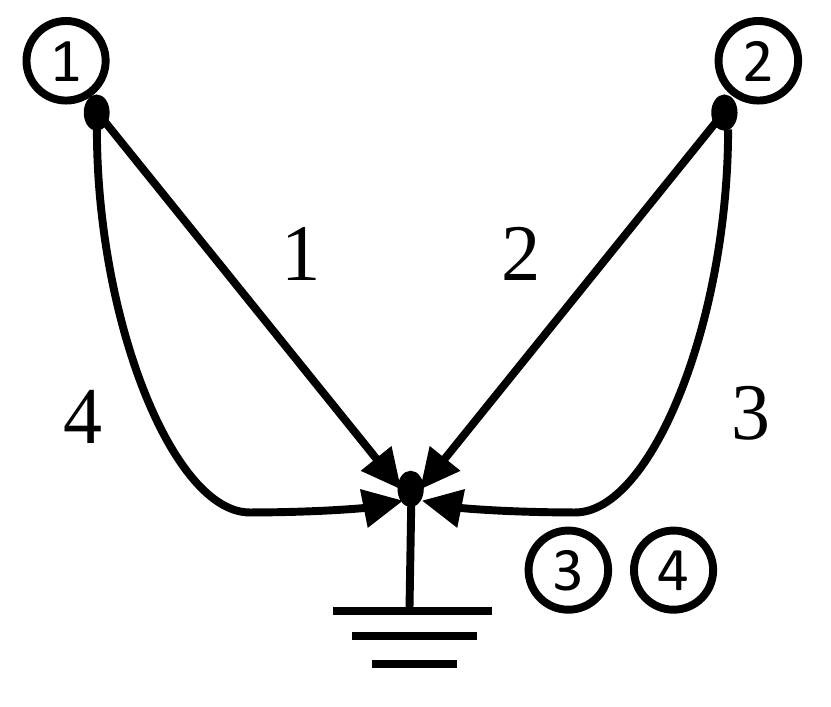}
		\label{fig:second_case}} 
	\caption{The linear element circuit and connected digraph}
	\label{fig:example}
\end{figure}

To illustrate, consider the simple circuit shown in Figure \ref{fig:example}. It consists of three elements: a voltage source, an ideal transformer described by $v_1 = \frac{n_1}{n_2}v_2$, $i_2 =- \frac{n_1}{n_2}i_1$, and a linear element described by $v_3=Zi_3$. 

Applying the preceding steps, we can construct KCL, KVL and element constitutive equations as following
\begin{align}
\label{eq:circuit_kcl}
\text{KCL:} && AI = 0 \Leftrightarrow
\underbrace{\begin{bmatrix}
1 & 0 & 0 & 1 \\
0 & 1 & 1 & 0 \\
\end{bmatrix}}_\textbf{A} 
\underbrace{\begin{bmatrix}
i_1 \\
i_2 \\
i_3 \\
i_4 \\
\end{bmatrix}}_\textbf{i} 
= \underbrace{\begin{bmatrix}
0 \\
0 \\
\end{bmatrix}}_\textbf{0} \\
\label{eq:circuit_kvl}
\text{KVL:} & & v - A^TV = \Leftrightarrow
\underbrace{\begin{bmatrix}
v_1 \\
v_2 \\
v_3 \\
v_4 \\
\end{bmatrix}}_\textbf{v} 
- \underbrace{\begin{bmatrix}
1 & 0 \\ 
0 & 1 \\
0 & 1 \\
1 & 0 \\
\end{bmatrix}}_\textbf{A$^T$} 
\underbrace{\begin{bmatrix}
V_1 \\
V_2 \\
\end{bmatrix}}_\textbf{V} 
= \underbrace{\begin{bmatrix}
0 \\
0 \\
0 \\
0 \\
\end{bmatrix}}_\textbf{0}
\end{align}
\begin{align}
\label{eq:circuit_branch}
\text{Linear Element Equations:} &  \\ \nonumber
\begin{rcases}
n_2v_1 - n_1v_2 = 0 \\ 
n_1i_1 + n_2i_2 = 0 \\
v_3 - Zi_3 = 0 \\
v_4 =  E \\
\end{rcases} & \Leftrightarrow
\underbrace{\begin{bmatrix}
n_2 & -n_1 & 0 & 0 \\
0 & 0 & 0 & 0 \\
0 & 0 & 1 & 0 \\
0 & 0 & 0 & 1 \\
\end{bmatrix}}_\textbf{F$_v$} 
\underbrace{\begin{bmatrix}
v_1 \\
v_2 \\
v_3 \\
v_4 \\
\end{bmatrix}}_\textbf{v} \\ \nonumber
 + & \underbrace{\begin{bmatrix}
0 & 0 & 0 & 0 \\
n_1 & n_2 & 0 & 0 \\
0 & 0 & -Z & 0 \\
0 & 0 & 0 & 0 \\
\end{bmatrix}}_\textbf{F$_i$}  
\underbrace{\begin{bmatrix}
i_1 \\
i_2 \\
i_3 \\
i_4 \\
\end{bmatrix}}_\textbf{i}  
= \underbrace{\begin{bmatrix}
0 \\
0 \\
0 \\
E \\
\end{bmatrix}}_\textbf{u$_s$} \nonumber 
\end{align}
												
As described above, incidence matrix $A$ and corresponding two constant matrices $F_v,F_i$ are immediately identifiable in the construction. 
												
Therefore, we observe that the circuit is linear $iff$ its branch equations can be written in the form of $F_v, F_i$, and time-invariant $iff$ both $F_v, F_i$ are constant with repect to time. As we will see in the next section, standard power system  elements within the transmission network are linear and time-invariant; however, we wil argue that for models common used in OPF, generators and loads may be represented as nonlinear current sources/sinks. \\ 
												
\textbf{DEFINITION 3.1.} (The Tableau Matrix) Since \eqref{eq:circuit_kcl}, \eqref{eq:circuit_kvl} and (\ref{eq:circuit_branch}) which constitute the tableau equation consist of a system of linear equations, it is convenient and more illuminating to rewrite them as a single matrix equation

\begin{align}
\label{eq:circuit_tableau}
\Scale[0.95]{\underbrace{\begin{bmatrix}[cc|cccc|cccc]
0 & 0 & 0 & 0 & 0 & 0 & 1 & 0 & 0 & 1 \\
0 & 0 & 0 & 0 & 0 & 0 & 0 & 1 & 1 & 0 \\ \hline
-1 & 0 & 1 & 0 & 0 & 0 & 0 & 0 & 0 & 0 \\  
0 & -1 & 0 & 1 & 0 & 0 & 0 & 0 & 0 & 0 \\
0 & -1 & 0 & 0 & 1 & 0 & 0 & 0 & 0 & 0 \\
-1 & 0 & 0 & 0 & 0 & 1 & 0 & 0 & 0 & 0 \\ \hline
0 & 0 & n_2 & -n_1 & 0 & 0 & 0 & 0 & 0 & 0 \\
0 & 0 & 0 & 0 & 0 & 0 & n_1 & n_2 & 0 & 0 \\
0 & 0 & 0 & 0 & 0 & 1 & 0 & 0 & -Z & 0 \\
0 & 0 & 0 & 0 & 0 & 0 & 0 & 0 & 0 & 0 \\
\end{bmatrix}}_\textbf{T}}
\Scale[0.95]{\underbrace{\begin{bmatrix}[c]
V_1 \\
V_2 \\  \hline
v_1 \\ 
v_2 \\
v_3 \\
v_4 \\ \hline
i_1 \\
i_2 \\
i_3 \\
i_4 \\
\end{bmatrix}}_\textbf{x}}
\Scale[0.95]{=\underbrace{\begin{bmatrix}[c]
0 \\
0 \\ \hline
0 \\
0 \\
0 \\
0 \\ \hline
0 \\
0 \\
0 \\
\Scale[0.9]{E} \\
\end{bmatrix}}_\textbf{u}}
\end{align}
															
The tableau matrix $T$ is as it is shown in (\ref{eq:circuit_tableau}) often very sparse, thereby allowing highly efficient numerical algorithms with a computer programming language. This is also why it is called sparse tableau analysis. In general, it can be recast into the compact matrix form with $\mathbf{0}$ and $\mathbf{I}$ denoting a zero and a unit matrix of appropriate dimension. 
\begin{align}
\label{eq:compact_tableau}
\text{Linear Tableau Formulation:} \nonumber &  \\
\underbrace{\begin{bmatrix}
\mathbf{0} & \mathbf{0} & A    \\
-A^T & \mathbf{I} & \mathbf{0} \\ 
\mathbf{0} & F_v & F_i   \\
\end{bmatrix}}_\textbf{T} &
\underbrace{\begin{bmatrix}
V \\
v \\ 
i \\
\end{bmatrix}}_\textbf{x} 
=\underbrace{\begin{bmatrix}
0 \\
0 \\
u_s \\
\end{bmatrix}}_\textbf{u}  
\end{align} 
																		
\section{Sparse Tableau Formulation for Power System Networks}
\label{chap:3.3}
This section carefully manipulates and applies the process of Sparse Tableau Formulation described in the previous section for power system networks. Standard power system network elements are represented as a two-port element and transmission matrix representation is employed to build branch equations for each network element. With all variables (port voltages $v$, port currents $i$, node voltages $V$ and node current $I$) defined, the corresponding incidence matrix $A$ is constructed to impose linear KCL and KVL equations. \\     

\textbf{PROPOSITION 1.} (Sparse Tableau Formulation for power system network) Standard power system networks can be cast into the compact matrix form using the sparse tableau matrix \eqref{eq:compact_tableau} with time-invariant circuit elements implying constant $F_v$, $F_i$ matrices, with each network element being a non-independent source and with each bus having current source implying $u_s = I$. This Sparse Tableau Formulation is shown as following \\

\noindent\text{Sparse Tableau Formulation for power system network:}
\begin{align}
\label{eq:compact_tableau_power}
\underbrace{\begin{bmatrix}
\mathbf{0} & \mathbf{0} & A    \\
-A^T & \mathbf{I} & \mathbf{0} \\ 
\mathbf{0} & F_v & F_i   \\
\end{bmatrix}}_\textbf{T} 
\underbrace{\begin{bmatrix}
V \\
v \\ 
i \\
\end{bmatrix}}_\textbf{x} 
=\underbrace{\begin{bmatrix}
I \\
0 \\
0 \\
\end{bmatrix}}_\textbf{u}  
\end{align}
where $I$ represents externally injected current source from generators or loads. Reasoning for the proof of this proposition is discussed systematically from the next section.
			
\subsection{Network Elements Modeling}
			
\begin{figure}[!htb]
\centering
\includegraphics[height=1in,width=2in]{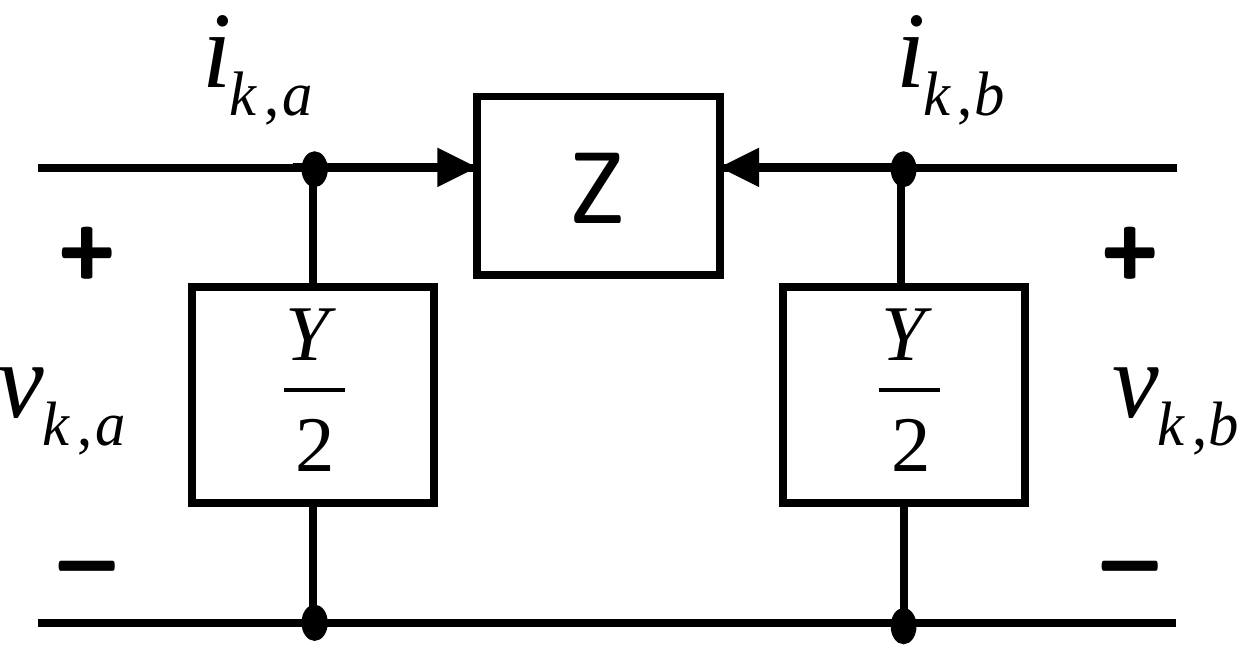}
\caption{$\pi$-equivalent circuit for transmission line}
\label{fig:line}
\end{figure}

$Proof.$ First step to prove the \textbf{PROPOSITION 1} is to consider the standard power system network circuit elements, which are a two-port element as depicted in Figure \ref{fig:twoport} to construct branch equations. Typical examples of two-port network elements would be transmisson lines and transformers. In OPF applications, transmission lines are often considered in terms of their $\pi$-equivalent circuit, rather than in the two-port [ABCD]-transmission matrix.   Hence, their data is typically provided as the three real-value parameters $R, X$ and $B$, with associated complex series impedance for line given by $Z = R+jX$ and shunt $Y = jB$.  At the sending end and receiving end ports, terminal behavior equivalent to the two-port may be captured in a circuit composed only of of simpler two-terminal elements \cite{BergenVittal2000}, as shown in Figure \ref{fig:line}. 

However, in standard power systems textbook presentations, one then recovers a two-port constitutive relation consistent wth \eqref{eq:general_element} by constructing the ``Transmission Matrix''.
    
\begin{align}
\label{eq:TM}
\begin{bmatrix} \nonumber
v_{k,a} \\
i_{k,a} \\ 
\end{bmatrix} & =
\begin{bmatrix}
1 + \frac{ZY}{2} & Z \\
Y(1 + \frac{ZY}{4}) & 1 + \frac{ZY}{2} \\ 
\end{bmatrix} 
\begin{bmatrix}
v_{k,b} \\
-i_{k,b} \\ 
\end{bmatrix} \\ 
& \Longrightarrow 
\underbrace{\begin{bmatrix}[cc|cc]
1 & 0 & -(1 + \frac{ZY}{2}) & Z \\ 
0 & 1 & -Y(1 + \frac{ZY}{4}) & (1 + \frac{ZY}{2}) \\
\end{bmatrix}}_\textbf{Representative circuit element $F_k$ $\in$ $\mathbb{C}^{2\times 4}$} 
\begin{bmatrix}
v_{k,a} \\
i_{k,a} \\ 
v_{k,b} \\
i_{k,b} \\
\end{bmatrix} =
\begin{bmatrix}
0 \\
0 \\ 
\end{bmatrix}
\end{align}

To adapt the equations (\ref{eq:TM}) to the Sparse Tableau Formulation (\ref{eq:compact_tableau_power}), we can re-write the branch equation as \\ 

\noindent \text{Linear Element Equation for transmission line:} 
\begin{align}
\label{eq:TM2}
& 
\begin{bmatrix}[cc|cc] \nonumber
1 & -(1 + \frac{ZY}{2}) & 0 & Z \\ 
0 & -Y(1 + \frac{ZY}{4}) & 1 & (1 + \frac{ZY}{2}) \\
\end{bmatrix} 
\begin{bmatrix}
v_{k,a} \\
v_{k,b} \\ 
i_{k,a} \\
i_{k,b} \\
\end{bmatrix} =
\begin{bmatrix}
0 \\
0 \\ 
\end{bmatrix} \\ 
\Scale[0.85]{\Longrightarrow} &
\underbrace{\begin{bmatrix}
1 & -(1 + \frac{ZY}{2}) \\ 
0 & -Y(1 + \frac{ZY}{4}) \\ 
\end{bmatrix}}_\textbf{F$_{k,v}$}
\begin{bmatrix}
v_{k,a} \\
v_{k,b} \\ 
\end{bmatrix} +
\underbrace{\begin{bmatrix}
0 & Z \\ 
1 & (1 + \frac{ZY}{2}) \\
\end{bmatrix}}_\textbf{F$_{k,i}$} 
\begin{bmatrix}
i_{k,a} \\
i_{k,b} \\ 
\end{bmatrix} =
\begin{bmatrix}
0 \\
0 \\ 
\end{bmatrix} 
\end{align}

thereby obtain the corresponding constant $F_v$ and $F_i$ matrices for the network element of transmission line. The other typical network element is a transformer. Here, voltage gain of transformer is expressed as complex scalar $T$ to account for phase shifting transformers (real-value voltage gain is for an ideal step up/down transformer allowing only voltage magnitude to change). Then, the corresponding transmission matrix representation is 
\begin{align}
\label{eq:TR}
\begin{bmatrix}
v_{k,a} \\
i_{k,a} \\ 
\end{bmatrix} & =
\begin{bmatrix}
T & 0 \\
0 & \frac{1}{T^*} \\ 
\end{bmatrix} 
\begin{bmatrix}
v_{k,b} \\
-i_{k,b} \\ 
\end{bmatrix} 
%& \Longrightarrow 
%\begin{bmatrix}[cc|cc]
%1 & 0 & -S & 0 \\ 
%0 & 1 & 0 & \frac{1}{S^*} \\
%\end{bmatrix}
%\begin{bmatrix}
%v_{k,a} \\
%i_{k,a} \\ 
%v_{k,b} \\
%i_{k,b} \\
%\end{bmatrix} =
%\begin{bmatrix}
%0 \\
%0 \\ 
%\end{bmatrix}
\end{align}
This can be equivalently re-written as \\

\noindent \text{Linear Element Equation for transformer:}
\begin{align}
\label{eq:TR2}
\underbrace{\begin{bmatrix}
1 & -T \\ 
0 & 0 \\ 
\end{bmatrix}}_\textbf{F$_{k,v}$} & 
\begin{bmatrix}
v_{k,a} \\
v_{k,b} \\ 
\end{bmatrix} +
\underbrace{\begin{bmatrix}
0 & 0 \\ 
1 & \frac{1}{T^*} \\
\end{bmatrix}}_\textbf{F$_{k,i}$} 
\begin{bmatrix}
i_{k,a} \\
i_{k,b} \\ 
\end{bmatrix} =
\begin{bmatrix}
0 \\
0 \\ 
\end{bmatrix} 
\end{align}

For each of these cases, it is straightforward to identify the corresponding constant $F_v$ and $F_i$ matrices. Note that for these very typical network elements, \eqref{eq:TM2} and \eqref{eq:TR2} contain only constant coefficients, have no independent sources, and therefore are simple linear functions.   

\subsection{Construction of the incidence matrix A}

Remaning constraints are simple linear expressions imposing KVL and KCL interconnection constraints. Since a node-to-element incident matrix $A$ is defined over all network elements, we need to organize all network element variables (port voltages and port currents):
\begin{align*}
v \triangleq \begin{bmatrix}[c]
v_{1,a} \\ 
v_{1,b} \\ \hline
\colon \\ 
\colon \\ \hline
v_{l,a} \\ 
v_{l,b}  \\ 
\end{bmatrix}, \quad
i \triangleq \begin{bmatrix}[c]
i_{1,a} \\ 
i_{1,b} \\ \hline
\colon \\ 
\colon \\ \hline
i_{l,a} \\ 
i_{l,b} \\ 
\end{bmatrix}
\end{align*}

Thus, $v,i \in \mathbb{C}^{2l}$ where $l$ is number of network elements. Goal of KCL is to efficiently to assemble the right hand side of the general current balance equation. To this end, the incidence matrix is then composed entirely of values of 1 or -1 or 0 by the following rule:
\begin{equation}
\label{eq:inc}
A(j,r) \in \mathbb{R}^{N \times 2l} \triangleq \begin{cases}
1, & \text{if $r$th component of $i$ corresponds to} \\
&  \text{an elements' sending or receiving} \\ 
& \text{terminal leaving node $j$} \\ \hline
-1, & \text{if $r$th component of $i$ corresponds to} \\
&  \text{an elements' sending or receiving} \\ 
& \text{terminal entering node $j$} \\ \hline
0, & \text{otherwise} 
\end{cases}
\end{equation}

Therefore, the current conservation law of KCL is written simply as 
\begin{align}
\label{eq:kcl}
I - Ai = 0 \,\,\, \in \mathbb{C}^{N}
\end{align}
where $I \in \mathbb{C}^{N}$ is the node complex current injection from generators or loads; $i \in \mathbb{C}^{2l}$ is the complex branch current carried away from node by network elements. We can also use $A$ to relate port voltages to bus voltages in a manner that guarantees KVL is automatically satisfied. Similarly, linear voltage law of KVL is written as  
\begin{align}
\label{eq:kvl}
v - A^TV = 0 \,\,\, \in \mathbb{C}^{2l}
\end{align}
where $V \in \mathbb{C}^{N}$ is bus voltages. The equation (\ref{eq:kvl}) is to assign the correct bus voltage to any port voltage of a port connected to that bus. Now, to construct the sparse tableau matrix \eqref{eq:compact_tableau_power}, $F_v$ and $F_i$ can be defined as
 
\begin{align}
F_v = \Scale[0.8]{\begin{bmatrix}
F_{1,v} & 0 & \cdots & 0 \\ 
0 & F_{2,v} & 0 & \colon \\  
\colon & \colon & \ddots & 0 \\ 
0 & \cdots & \cdots & F_{l,v}\\ 
\end{bmatrix}} \in \Scale[0.9]{\mathbb{C}^{2l\times 2l}}  \,\,
F_i = \Scale[0.8]{\begin{bmatrix}
F_{1,i} & 0 & \cdots & 0 \\ 
0 & F_{2,i} & 0 & \colon \\  
\colon & \colon & \ddots & 0 \\ 
0 & \cdots & \cdots & F_{l,i}\\ 
\end{bmatrix}} \in \Scale[0.9]{\mathbb{C}^{2l\times 2l}}
\end{align} 
where $F_v$ and $F_i$ are block diagonal matrix composed of previously described $F_{k,v}$, $F_{k,i}$ matrices. Finally, with marices $F_v$, $F_i$, $A$ and variables $v$, $i$, $V$, $I$ as defined above, we can describe power system networks using the Sparse Tableau Formulation by \\

\noindent\text{Sparse Tableau Formulation for power system:}
\begin{align}
\label{eq:final_tableau_power}
\underbrace{\begin{bmatrix}
\mathbf{0} & \mathbf{0} & A    \\
-A^T & \mathbf{I} & \mathbf{0} \\ 
\mathbf{0} & F_v & F_i   \\
\end{bmatrix}}_\textbf{T} 
\underbrace{\begin{bmatrix}
V \\
v \\ 
i \\
\end{bmatrix}}_\textbf{x} 
=\underbrace{\begin{bmatrix}
I \\
0 \\
0 \\
\end{bmatrix}}_\textbf{u} \qquad\IEEEQEDhere
\end{align} 

\begin{figure}[!htb]
	\centering
	\includegraphics[height=1in,width=1in]{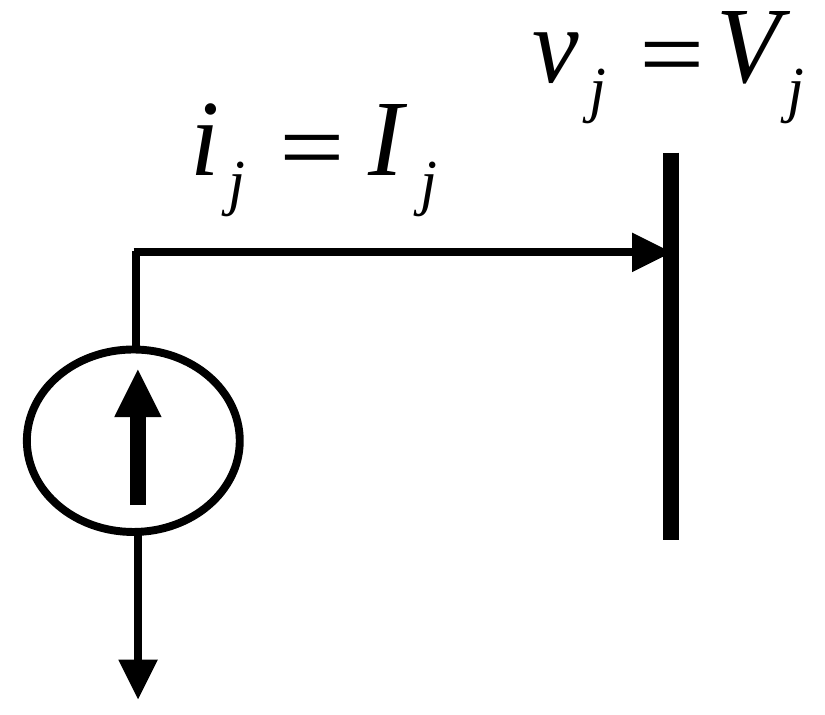}
	\caption{Nonlinear current source element as one port element}
	\label{fig:oneport}
\end{figure}

Notice that here we define current source elements $I$ and this introduces a special class of nonlinear one-port element in Figure \ref{fig:oneport} to describe power system network for power flow analysis. 
Then, nonlinear element equation as a equation (\ref{eq:general_element}) for current source $I_j$ for bus $j$ can be defined by 
\begin{align}
\text{Nonlinear Element Equations:} & \nonumber \\ 
f_j (v_{j}, i_{j}) = 0 \triangleq i_j = & \frac{(S_{g,j} - S_{d,j})^*}{v_j^*} 
\label{eq:nonlinear_ele}
\end{align}
\text{Notice that $S_j = S_{g,j} - S_{d,j}$, $i_j = I_j$, and $v_j = V_j$} implying
\begin{align}
\Longrightarrow & I_j - \frac{S_j^*}{V_j^*} = 0 
\end{align}
where $S_{g,j}$ and $S_{d,j}$ are specified apparent power generation and load at bus $j$. Notice that equation (\ref{eq:nonlinear_ele}) is similar to approach in traditional Gauss-Seidel formulation of power flow \cite{BergenVittal2000} and is equivalent to $S_j = V_jI_j^*$, which is typical ``power balance equation''. %In this paper, we stick to the form (\ref{eq:nonlinear_ele}) since it clarifies our concept of Sparse Tableau Formulation and port representation of network elements. 

\subsection{Illustrative example with three-bus system}

\begin{figure}[!htb]
	\centering
	\includegraphics[height=2in,width=3in]{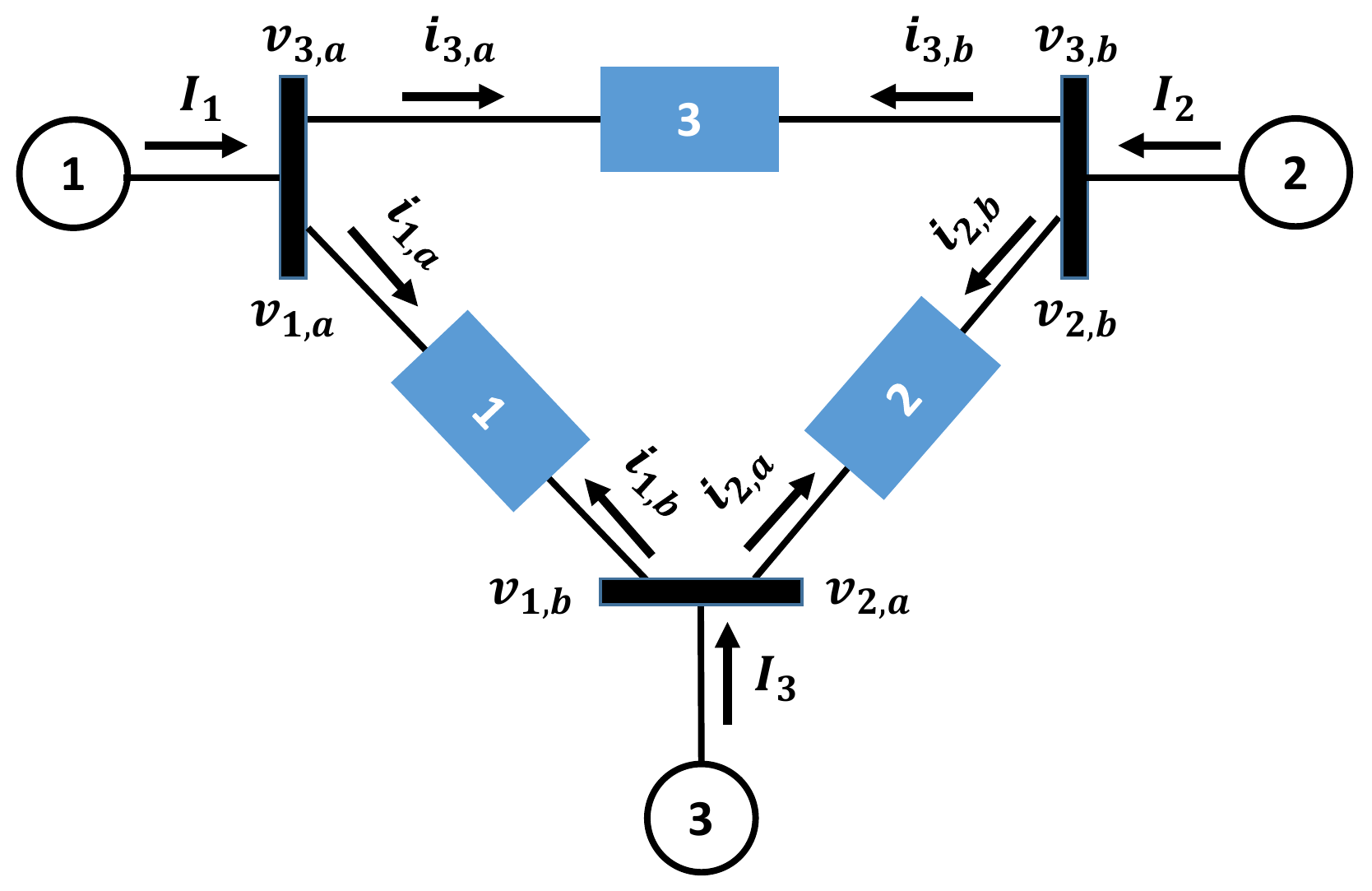}
	\caption{Three-bus system for sparse tableau formulation}
	\label{fig:three}
\end{figure}

This section details the STF for power system networks by providing an illustrative example with three-bus system depicted in Figure \ref{fig:three}. 
Notice that this system consists of three buses as \circled{1}, \circled{2}, \circled{3} and three network elements, three transmisson lines as \boxed{1}, \boxed{2}, \boxed{3} with \boxed{1} : line from bus \circled{1} to \circled{3}, \boxed{2} : line from bus \circled{3} to \circled{2} and \boxed{3} : line from bus \circled{1} to \circled{2}. Based on the rule (\ref{eq:inc}), we can construct the incident matrix $A$ as

\begin{align}
\label{eq:exam_inc}
A =
\Scale[0.9]{\begin{blockarray}{cc|cc|ccc}
\boxed{1},a & \boxed{1},b & \boxed{2},a & \boxed{2},b & \boxed{3},a & \boxed{3},b \\
\begin{block}{(cc|cc|cc)c} 
1 & 0 & 0 & 0 & 1 & 0 & \circled{1} \\
0 & 0 & 0 & 1 & 0 & 1 & \circled{2} \\
0 & 1 & 1 & 0 & 0 & 0 & \circled{3} \\
\end{block}
\end{blockarray}}
\end{align}

Notice that each row corresponds to the node and each two-column corresponds to port $a$ or $b$ of the network element. Then, branch equations for each network element can be constructed with $F_{\scalebox{0.6}{$\boxed{1}$},v}$, $F_{\scalebox{0.6}{$\boxed{2}$},v}$, and $F_{\scalebox{0.6}{$\boxed{3}$},v}$ 
\begin{align}
\label{eq:exam_FV}
F_{\scalebox{0.6}{$\boxed{1}$},v} = 
\Scale[0.9]{\begin{bmatrix}
1 & -(1 + \frac{Z_{\scalebox{0.6}{$\boxed{1}$}}Y_{\scalebox{0.6}{$\boxed{1}$}}}{2}) \\ 
0 & -Y_{\scalebox{0.6}{$\boxed{1}$}}(1 + \frac{Z_{\scalebox{0.6}{$\boxed{1}$}}Y_{\scalebox{0.6}{$\boxed{1}$}}}{4}) \\ 
\end{bmatrix}} \\
F_{\scalebox{0.6}{$\boxed{2}$},v} = 
\Scale[0.9]{\begin{bmatrix}	
1 & -(1 + \frac{Z_{\scalebox{0.6}{$\boxed{2}$}}Y_{\scalebox{0.6}{$\boxed{2}$}}}{2}) \\ 
0 & -Y_{\scalebox{0.6}{$\boxed{2}$}}(1 + \frac{Z_{\scalebox{0.6}{$\boxed{2}$}}Y_{\scalebox{0.6}{$\boxed{2}$}}}{4}) \\ 
\end{bmatrix}}  \\ 
F_{\scalebox{0.6}{$\boxed{3}$},v} = 
\Scale[0.9]{\begin{bmatrix}	
1 & -(1 + \frac{Z_{\scalebox{0.6}{$\boxed{3}$}}Y_{\scalebox{0.6}{$\boxed{3}$}}}{2}) \\ 
0 & -Y_{\scalebox{0.6}{$\boxed{3}$}}(1 + \frac{Z_{\scalebox{0.6}{$\boxed{3}$}}Y_{\scalebox{0.6}{$\boxed{3}$}}}{4})\\ 
\end{bmatrix}}
\end{align}
Similarly, $F_{\scalebox{0.6}{$\boxed{1}$},i}$, $F_{\scalebox{0.6}{$\boxed{2}$},i}$, and $F_{\scalebox{0.6}{$\boxed{3}$},i}$ can be constructed as
\begin{align}
\label{eq:exam_FI}
F_{\scalebox{0.6}{$\boxed{1}$},i} = 
\Scale[0.9]{\begin{bmatrix}
0 & Z_{\scalebox{0.6}{$\boxed{1}$}} \\ 
1 & (1 + \frac{Z_{\scalebox{0.6}{$\boxed{1}$}}Y_{\scalebox{0.6}{$\boxed{1}$}}}{2}) \\ 
\end{bmatrix}} \\
F_{\scalebox{0.6}{$\boxed{2}$},i} = 
\Scale[0.9]{\begin{bmatrix}	
0 & Z_{\scalebox{0.6}{$\boxed{2}$}} \\ 
1 & (1 + \frac{Z_{\scalebox{0.6}{$\boxed{2}$}}Y_{\scalebox{0.6}{$\boxed{2}$}}}{2}) \\ 
\end{bmatrix}}  \\
F_{\scalebox{0.6}{$\boxed{3}$},i} = 
\Scale[0.9]{\begin{bmatrix}	
0 & Z_{\scalebox{0.6}{$\boxed{3}$}} \\ 
1 & (1 + \frac{Z_{\scalebox{0.6}{$\boxed{3}$}}Y_{\scalebox{0.6}{$\boxed{3}$}}}{2}) \\ 
\end{bmatrix}}
\end{align}
where $Z_{\scalebox{0.6}{$\boxed{k}$}}$, $Y_{\scalebox{0.6}{$\boxed{k}$}}$ are parameters for the network element $\boxed{k}$. Here, we only consider transmission lines for the network element, but $F_v$ and $F_i$ for transformers can be easily constructed by (\ref{eq:TR2}). Next, we can define all corresponding variables as 
\begin{align*}
v \triangleq \begin{bmatrix}
v_{1,a} \\ 
v_{1,b} \\  
v_{2,a} \\ 
v_{2,b} \\ 
v_{3,a} \\ 
v_{3,b} \\ 
\end{bmatrix}, \quad
i \triangleq \begin{bmatrix}
i_{1,a} \\ 
i_{1,b} \\  
i_{2,a} \\ 
i_{2,b} \\ 
i_{3,a} \\ 
i_{3,b} \\ 
\end{bmatrix}, \quad
V \triangleq \begin{bmatrix}
V_{1} \\ 
V_{2} \\ 
V_{3} \\ 
\end{bmatrix}, \quad
I \triangleq \begin{bmatrix}
I_{1} \\ 
I_{2} \\  
I_{3} \\ 
\end{bmatrix}
\end{align*}
\def\block(#1,#2)#3{\multicolumn{#2}{c}{\multirow{#1}{*}{$#3$}}}
\def\bblock(#1,#2)#3{\multicolumn{#2}{c|}{\multirow{#1}{*}{$#3$}}}
and then the STF for this three-bus system network is expressed as  
\begin{align}
\label{eq:exam_tableau}
\Scale[0.95]{\underbrace{\begin{bmatrix}[ccc|cccccc|cccccc]
0 & 0 & 0 & 0 & 0 & 0 & 0 & 0 & 0 & 1 & 0 & 0 & 0 & 1 & 0 \\
0 & 0 & 0 & 0 & 0 & 0 & 0 & 0 & 0 & 0 & 0 & 0 & 1 & 0 & 1 \\ 
0 & 0 & 0 & 0 & 0 & 0 & 0 & 0 & 0 & 0 & 1 & 1 & 0 & 0 & 0 \\ \hline 
-1 & 0 & 0 & 1 & 0 & 0 & 0 & 0 & 0 & 0 & 0 & 0 & 0 & 0 & 0 \\
0 & 0 & -1 & 0 & 1 & 0 & 0 & 0 & 0 & 0 & 0 & 0 & 0 & 0 & 0 \\
0 & 0 & -1 & 0 & 0 & 1 & 0 & 0 & 0 & 0 & 0 & 0 & 0 & 0 & 0 \\ 
0 & -1 & 0 & 0 & 0 & 0 & 1 & 0 & 0 & 0 & 0 & 0 & 0 & 0 & 0 \\
-1 & 0 & 0 & 0 & 0 & 0 & 0 & 1 & 0 & 0 & 0 & 0 & 0 & 0 & 0 \\
0 & -1 & 0 & 0 & 0 & 0 & 0 & 0 & 1 & 0 & 0 & 0 & 0 & 0 & 0 \\ \hline
0 & 0 & 0 & \block(2,2){F_{\scalebox{0.6}{$\boxed{1}$},v}} & 0 & 0 & 0 & 0 & \block(2,2){F_{\scalebox{0.6}{$\boxed{1}$},i}} & 0 & 0 & 0 & 0 \\
0 & 0 & 0 &  &  & 0 & 0 & 0 & 0 &  &  & 0 & 0 & 0 & 0\\
0 & 0 & 0 & 0 & 0 & \block(2,2){F_{\scalebox{0.6}{$\boxed{2}$},v}} & 0 & 0 & 0 & 0 & \block(2,2){F_{\scalebox{0.6}{$\boxed{2}$},i}} & 0 & 0\\
0 & 0 & 0 & 0 & 0 &  &  & 0 & 0 & 0 & 0 &  &  & 0 & 0\\
0 & 0 & 0 & 0 & 0 & 0 & 0 & \bblock(2,2){F_{\scalebox{0.6}{$\boxed{3}$},v}} & 0 & 0 & 0 & 0 & \block(2,2){F_{\scalebox{0.6}{$\boxed{3}$},i}} \\
0 & 0 & 0 & 0 & 0 & 0 & 0 &  &  & 0 & 0 & 0 & 0 &  &  \\
\end{bmatrix}}_\textbf{T}} 
\Scale[0.95]{\underbrace{\begin{bmatrix}[c]
V_1 \\
V_2 \\ 
V_3 \\ \hline
v_{1,a} \\
v_{1,b} \\
v_{2,a} \\
v_{2,b} \\ 
v_{3,a} \\
v_{3,b} \\ \hline
i_{1,a} \\
i_{1,b} \\
i_{2,a} \\
i_{2,b} \\ 
i_{3,a} \\
i_{3,b} \\ 
\end{bmatrix}}_\textbf{x}} 
\Scale[0.95]{=\underbrace{\begin{bmatrix}[c]
I_1 \\
I_2 \\ 
I_3 \\ \hline
0 \\
0 \\
0 \\
0 \\ 
0 \\ 
0 \\ \hline
0 \\
0 \\
0 \\
0 \\
0 \\
0 \\
\end{bmatrix}}_\textbf{u}}  
\end{align}

\section{Property of Sparse Tableau Formulation}
\label{chap:3.4}

\subsection{Relationship with $\textbf{Y}_{bus}$}
Standard power system network uses the admittance matrix $Y_{bus}$ to model power system networks. $Y_{bus} \in \mathbb{C}^{N\times N} $ represents the nodal admittance of the buses in power system networks by writing KCL equations at nodes/buses in terms of node voltages, which is referred to ``nodal analysis''. A more general prodecure for constructing the $Y_{bus}$ can be found \cite{SadikuAlexander2011}, \cite{Grainger1994}. Here, we analyze the relationship between the $Y_{bus}$ and STF by constructing the $Y_{bus}$ from the STF.\\

\textbf{LEMMA 1.} (Y$_{bus}$ and Sparse Tableau Formulation) Standard modeling via $Y_{bus}$, used in most OPF formulations, represents algebraic elimination of variables from the Sparse Tableau Fomulation with a restrictive assumption.\\

\noindent $Proof.$ Observe that the $Y_{bus}$ describes a linear relation of bus voltage to bus current as:
\begin{align}
\label{eq:Ybus2}
I = Y_{bus}V
\end{align}
To relate the $Y_{bus}$ to the Sparse Tableau Formulation,  we perform an algebraic reduction to obtain \eqref{eq:Ybus2} from the Sparse Tableau Formulation. Consider the full Sparse Tableau Formulation:   
\begin{align} 
\begin{cases}
Ai = I \\
v = A^TV \\ 
F_v v + F_i i = 0 \\
\end{cases}   
\end{align} 
With these equations, we can solve for $i$ such that 
\begin{align}
i = -(F_i)^{-1} \cdot F_v \cdot v
\end{align} 
Then, substitute the KVL equation $v=A^TV$, yielding 
\begin{align}
i & = -(F_i)^{-1} \cdot F_v \cdot A^TV \quad \text{by KCL $I=Ai$} \\
I & = \underbrace{-A\cdot (F_i)^{-1} \cdot F_v \cdot A^T}_\textbf{Y$_{bus}$}V
\end{align} 
This demonstrates that the $Y_{bus}$ formulation may be viewed as a simple algebraic elimination of variables from the STF and this elimination requires a restrictive assumption that $F_i$ is a full rank matrix, which also requires a full rank condition of each $F_{k,i}$ matrix. $\IEEEQEDhere$

This lemma shows that the STF is more general modeling approach, in the sense that $Y_{bus}$ is obtained from the STF only with the restrictive assumption that $F_i$ is invertible. Important network elements such as an ideal transformer and circuit breaker fail this requirement.   

\subsection{Non-typical Network Elements using the Sparse Tableau Formulation}

Part of the argument for the Sparse Tableau Formulation lies in the convenience with which non-typical network elements may be treated, in contrast to the $Y_{bus}$ formulation, that often requires ``tricks" in handling such elements. In this section, we provide three examples for non-typical network elements: 1) Ideal Transformer, 2) Circuit Breaker, and 3) Three-Winding Transformer. \\

\textbf{Non-typical Element 1.} (Ideal Transformer) Consider an ideal transformer model \eqref{eq:TR2}
\begin{align}
\text{Ideal transformer:} \nonumber & \\ 
\underbrace{\begin{bmatrix}
1 & -T \\ 
0 & 0 \\ 
\end{bmatrix}}_\textbf{F$_{v}$} &
\begin{bmatrix}
v_a \\
v_b \\ 
\end{bmatrix} +
\underbrace{\begin{bmatrix}
0 & 0 \\ 
1 & \frac{1}{T^*} \\
\end{bmatrix}}_\textbf{F$_{i}$} 
\begin{bmatrix}
i_a \\
i_b \\ 
\end{bmatrix} =
\begin{bmatrix}
0 \\
0 \\ 
\end{bmatrix} 
\end{align}

It can be easily checked that $F_i$ is not invertible, which implies that modeling approach with the $Y_{bus}$ fails to model the ideal transformer as a stand-alone element. Power system specialists will recognize that in $Y_{bus}$ formulations, and ideal transformer is never treated ``stand-alone." Rather, series or shunt impedances representing leakage and/or magnetizing reactances are always added.  \\

\textbf{Non-typical Element 2.} (Circuit Breaker) More importantly, consider the circuit representation of a circuit breaker with binary integer parameter $\gamma$ indicating switch position.
\begin{align}
& \text{Circuit breaker closed, $\gamma = 1$ :} \nonumber 
\begin{cases}
v_a - v_b = 0 \\
i_a + i_b = 0 \\ 
\end{cases}  \\
& \Longrightarrow  
\underbrace{\begin{bmatrix}
1 & -1 \\ 
0 & 0 \\ 
\end{bmatrix}}_\textbf{F$_{v}$}
\begin{bmatrix}
v_a \\
v_b \\ 
\end{bmatrix} +
\underbrace{\begin{bmatrix}
0 & 0 \\ 
1 & 1 \\
\end{bmatrix}}_\textbf{F$_{i}$} 
\begin{bmatrix}
i_a \\
i_b \\ 
\end{bmatrix} =
\begin{bmatrix}
0 \\
0 \\ 
\end{bmatrix} \\
& \text{Circuit breaker open, $\gamma = 0$ :} \nonumber 
\begin{cases}
i_a = 0 \\
i_b = 0 \\ 
\end{cases} \\ 
& \Longrightarrow 
\underbrace{\begin{bmatrix}
0 & 0 \\ 
0 & 0 \\ 
\end{bmatrix}}_\textbf{F$_{v}$}
\begin{bmatrix}
v_a \\
v_b \\ 
\end{bmatrix} +
\underbrace{\begin{bmatrix}
1 & 0 \\ 
0 & 1 \\
\end{bmatrix}}_\textbf{F$_{i}$} 
\begin{bmatrix}
i_a \\
i_b \\ 
\end{bmatrix} =
\begin{bmatrix}
0 \\
0 \\ 
\end{bmatrix}  \\
%\end{align}
%\begin{align}
& \text{Circuit breaker :}  \nonumber \\
& \underbrace{\begin{bmatrix}
\gamma & -\gamma \\ 
0 & 0 \\ 
\end{bmatrix}}_\textbf{F$_{v}$}
\begin{bmatrix}
v_a \\
v_b \\ 
\end{bmatrix} +
\underbrace{\begin{bmatrix}
(1-\gamma) & 0 \\ 
\gamma & 1 \\
\end{bmatrix}}_\textbf{F$_{i}$} 
\begin{bmatrix}
i_a \\
i_b \\ 
\end{bmatrix} =
\begin{bmatrix}
0 \\
0 \\ 
\end{bmatrix}  
\end{align} 
Similarly, $F_i$ is not guaranteed invertible for both breaker positions, which implies that modeling approach with a fixed $Y_{bus}$ fails to model the circuit breaker. Instead, the $Y_{bus}$ based analysis must rely on ``topology processing," which may be viewed as a means to compute and switch between families of different $Y_{bus}$ matrices, depending on breaker settings. \\

\begin{figure}[htb!]
	\centering
	\subfigure[Structure of a three-winding transformer]
	{\includegraphics[height=1in, width=2.5in]{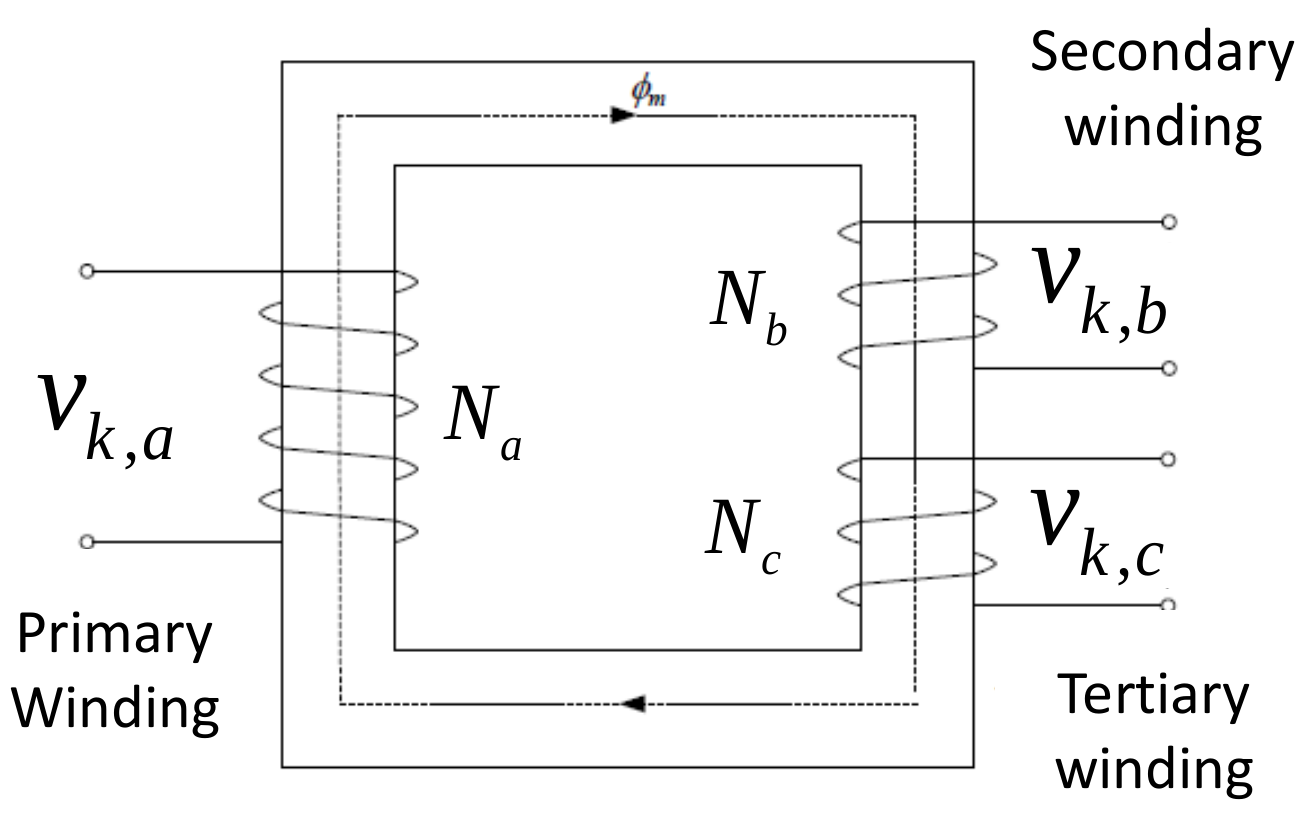}
		\label{fig:three_tr1}}
	\hfil
	\subfigure[Three-port representation]
	{\includegraphics[height=0.9in, width=2.5in]{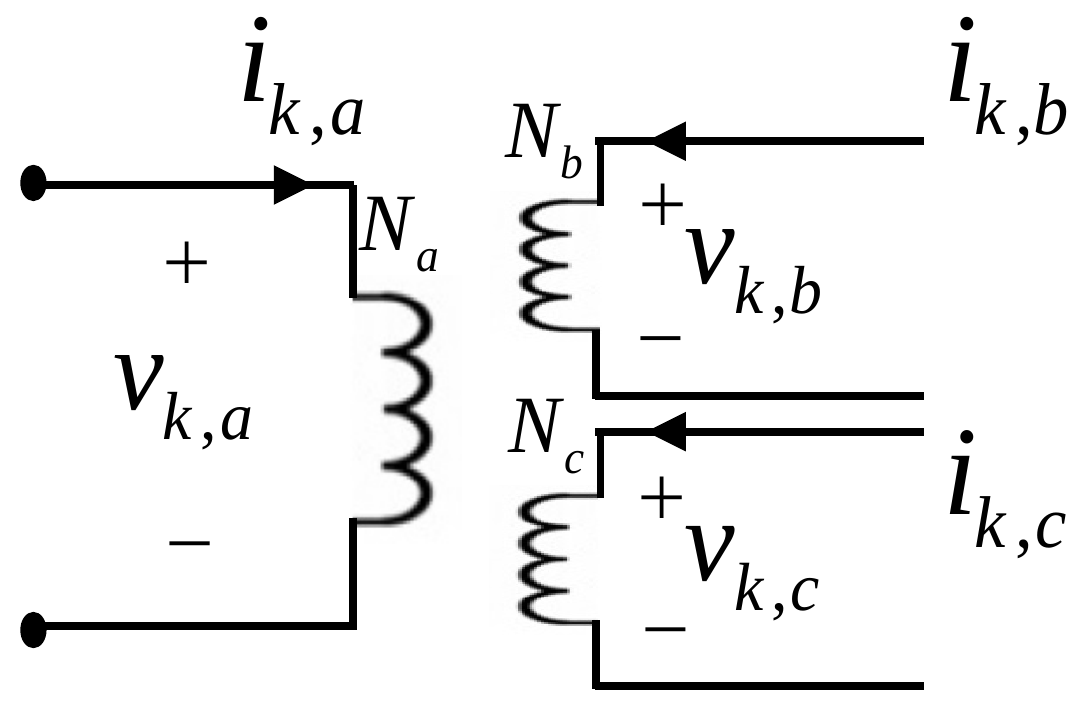}
		\label{fig:three_tr2}} 
	\caption{Schematic diagram of a transformer with three windings}
	\label{fig:three_tr}
\end{figure}

\textbf{Non-typical Element 3.} (Three-Winding Transformer) In power system networks, there is a small number of three-winding transformers best represented as three-ports depicted in Figure \ref{fig:three_tr}. Traditional $Y_{bus}$ approach to model the three-winding is to use the equivalent impedances $Z_p$, $Z_s$ and $Z_t$ \cite{Weedy1987}, \cite{OommenKohler1993}, \cite{Shaarbafi2014}, and the magnetizing current and core losses are modelled as a magnetizing branch. Arbitrary location of the magnetizing branch is one of the drawbacks of this approach. In addition, this modeling approach with the $Y_{bus}$ can create numeric difficulties, because in most large transformers the value of $Z_s$ is very small or even negative. The use of non-physical negative impedances, while common practice, is arguably undesirable and unrealistic.   

However, the modeling of three-winding transformer becomes simple and accurate with the Sparse Tableau Formulation. Consider the ideal conditions for a three-winding transformer (which are similar to those of two-winding ideal transformer):
\begin{align} 
& \Longrightarrow \begin{cases}
\frac{1}{N_a}v_{a} - \frac{1}{N_b}v_{b} = 0 \\
\frac{1}{N_b}v_{b} - \frac{1}{N_c}v_{c} = 0 \\ 
N_ai_{a} + N_bi_{b} + N_ci_{c} = 0 \\
\end{cases}   \nonumber
\end{align}
\begin{align}
\text{Linear Element Equations for three-winding transformer:} & \nonumber \\
\underbrace{\begin{bmatrix}
\frac{1}{N_a} & -\frac{1}{N_b} & 0 \\ 
0 & \frac{1}{N_b} & -\frac{1}{N_c} \\ 
0 & 0 & 0 \\ 
\end{bmatrix}}_\textbf{F$_{v}$}
\begin{bmatrix}
v_a \\
v_b \\
v_c \\  
\end{bmatrix}  + 
\underbrace{\begin{bmatrix}
0 & 0 & 0 \\ 
0 & 0 & 0 \\
N_a & N_b & N_c \\
\end{bmatrix}}_\textbf{F$_{i}$} 
\begin{bmatrix}
i_a \\
i_b \\ 
i_c \\ 
\end{bmatrix} =
\begin{bmatrix}
0 \\
0 \\
0 \\ 
\end{bmatrix} 
\end{align}
Thus, we can easily model a three-winding transformer as a three-port network element with the Sparse Tableau Formulation. 

\section{OPF Numeric Case Studies Employing the Sparse Tableau Formulation}
\label{chap:3.5}
This section applies the STF to optimal power flow problem. The optimal power flow determines decision variable values that produce an optimal operation point for the power system network in terms of a specified objective function with both network constraints and engineering constraints. The sum of individual generator cost functions is typically chosen as the objective function. Define an individual generator cost function as follow:
\begin{equation}
\tilde{c}_j(P_{g,j}) = \alpha (P_{g,j})^{2} + \beta P_{g,j} + \gamma,
\quad \forall j \in \mathbf{G} \\ 
\end{equation}

The OPF problem is non-convex due to the nature of the power flow equation \cite{LesieutreHiskens2005} and it is, in general, hard to obtain a guaranteed global solution. Non-convexity of the OPF problem has made solution techniques an ongoing research topic since the problem was first introduced by Carpentier in 1962 \cite{Carpentier1962}.  

\subsection{Optimal Power Flow with the Sparse Tableau Formulation}

Given a connected power system with $\mathbf{N}$ the set of all buses, $\mathbf{G}$ the set of all generators and $\mathbf{L}$ set of all transmission line, suppose that generators and loads connected to the network are specified by the active and reactive power they inject and withdraw at node respectively. So associated with each bus $j$, net complex power injection is $S_{j}$ = $P_j + jQ_j$ where $P_j = P_{g,j} - P_{d,j}$ and $Q_j = Q_{g,j} - Q_{d,j}$. In OPF, $P_j$ and $Q_j$ will be a decision variables for bus $j$ corresponding to a generator, and net complex power injection $S_j$ is given as nonlinear element equations with the equation (\ref{eq:nonlinear_ele}). Then, a representative OPF problem might take the following form: 
\begin{subequations}
\label{eq:sparse_opf}
\begin{equation}
\min_{\substack{P, Q, v, i, V, I}} \quad \quad \sum_{j \in \mathbf{G}} \tilde{c}_j\left(P_{g,j}\right)  \quad \text{subject to} \nonumber
\end{equation}
\begin{align}
\label{eq:sparse_branch_1}
\text{Linear Element:} \quad &  F_{v}v + F_{i}i = 0 \\
\text{KCL:} \quad &  I - Ai = 0 \\
\text{KVL:} \quad &  v - A^TV = 0 \\
\text{Nonlinear Element:} \quad &  S - V\odot (I)^* = 0 \\   
\text{Gen. Limit:} \quad &  P_{j}^{min} \leq P_{g,j} \leq P_{j}^{max} \\
\quad & Q_{j}^{min} \leq Q_{g,j} \leq Q_{j}^{max}, \, \forall j \in \mathbf{G} \\
\text{Vol. Limit:} \quad &  V_{j}^{min} \leq \lvert V_j \rvert \leq V_{j}^{max}, \, \forall j \in \mathbf{N} \\
\text{Line Limit:} \quad &  \lvert i_{k,a/b}\rvert \leq i_{k}^{max}, \,  \forall k \in \mathbf{L}
\end{align}
\end{subequations}

One of the benefit to use the STF for OPF problem is an explicit appearance of variable $i$, current flow on branch, so that thermal line limits with current magnitude as a superior choice \cite{CongRegulskiWallEtAl2015} can be easily incorporated without an additional effort to express the equation of current flows. \\

\textbf{REMARK 1.} (Experiment for Computational Time) Sparse Tableau Formulation of OPF problem is compared with two standard OPF problems; the polar power-voltage and rectangular current-voltage formulation selected as preffered formulations in terms of computational time \cite{ParkTangFerrisEtAl2016}. Problems are formulated in GAMS and solver KNITRO is selected based on the experience that there are a couple of advantages that KNITRO has over other solvers \cite{MurliLeonePardalosEtAl1998}.

\begin{table}[htb]
\renewcommand{\arraystretch}{1.3}
\centering
\resizebox{0.5\textwidth}{!}{%
\begin{tabular}{|c|c|c|c|c|c|c|}
\hline
\multicolumn{1}{|c|}{\multirow{2}{*}{}} & \multicolumn{2}{c|}{POLAR-$Y_{bus}$} & \multicolumn{2}{c|}{REC-IV-$Y_{bus}$} & \multicolumn{2}{c|}{REC-STF} \\ \cline{2-7} 
\multicolumn{1}{|c|}{} & Obj & Time & Obj & Time & Obj & Time \\ \hline\hline
\multirow{1}{*}{case118} & 129660.68 & 0.3sec & 129660.68 & 0.3sec & 129660.68 & 0.3sec \\ \hline
\multirow{1}{*}{case300}  & 719725.07 & 0.6sec & 719725.07 & 2sec & 719725.07 & 0.6sec \\ \hline 
\multirow{1}{*}{case2383wp} & 1868511.82 & 5.6sec & 1862367.02 & 5.2sec & 1862367.02 & 6.3sec \\ \hline
\multirow{1}{*}{case3012wp} & 2591706.57 & 5.8sec & 2582670.47 & 5.9sec & 2582670.47 & 6sec \\ \hline		
\multirow{1}{*}{case3120sp}  & 2142703.76 & 5.8sec & 2141532.10 & 6.3sec & 2141532.10 & 6.5sec \\ \hline 
\multirow{1}{*}{case3375wp} & 7412030.67 & 54sec & 7404635.99 & 11.7sec & 7404637.15 & 11.4sec \\ \hline
\noalign{\vskip 2mm} 
\end{tabular}
}\caption{Comparison of ACOPF problems}
\label{tab:sparse_acopf}
\end{table}

TABLE \ref{tab:sparse_acopf} shows the objective value and computational time of three ACOPF formulations. Notice that because REC-IV-$Y_{bus}$ and REC-STF limit current magnitudes rather than apparent power on lines, the solutions tend to be slightly different than the POLAR-$Y_{bus}$ formulation. Three formulations show that a similar computational time is required to solve the problem for most test cases, and POLAR-$Y_{bus}$ or REC-IV-$Y_{bus}$ or REC-STF is more superior than others for some cases, which demonstrates that the Sparse Tableau Formulation of OPF preserves the computational efficiency.

\section{Conclusion}
In this paper, we have addressed the Sparse Tableau Formulation (STF) for power system networks. After deriving the STF for power system networks, an illustrative example with three-bus system is given for developing of a concrete understanding. Then, relationship with the $Y_{bus}$ formulation is studied, and it has shown that the $Y_{bus}$ formulation is one of the special modeling approach from the STF with the restrictive assumption that $F_i$ matrix needs to be invertible. Specific network elements in which the STF becomes useful include three non-typical network elements: ideal transformer, circuit breaker and three-winding transformer. 

For computational comparison, optimal power flow problem is constructed with the STF. It shows the benefit of having current flow $i$ as a variables to impose thermal line limits and is compared with standard ACOPF formulation. The experiment demonstrates that the equivalent computational efficiency of the Sparse Tableau Formulation for OPF problem.

% or
								
								
%\appendices
%% you can choose not to have a title for an appendix
%% if you want by leaving the argument blank
%\section{}

\section*{Acknowledgment}
Work described here was supported by the Advanced Research Projects Agency-Energy (ARPA-E), U.S. Department of Energy, under Award Number DEAR0000717. The authors gratefully acknowledge this support; however, views and opinions of the author expressed herein do not necessarily state of reflect those of the United States Government or any agency thereof.
								
\ifCLASSOPTIONcaptionsoff
\newpage
\fi
								
								
								
%\bibliographystyle{IEEEtran}
								
%\begin{thebibliography}{1}
%\end{thebibliography}
								
\bibliographystyle{IEEEtran}
\bibliography{sparse}

\begin{thebibliography}{10}
\providecommand{\url}[1]{#1}
\csname url@samestyle\endcsname
\providecommand{\newblock}{\relax}
\providecommand{\bibinfo}[2]{#2}
\providecommand{\BIBentrySTDinterwordspacing}{\spaceskip=0pt\relax}
\providecommand{\BIBentryALTinterwordstretchfactor}{4}
\providecommand{\BIBentryALTinterwordspacing}{\spaceskip=\fontdimen2\font plus
\BIBentryALTinterwordstretchfactor\fontdimen3\font minus
  \fontdimen4\font\relax}
\providecommand{\BIBforeignlanguage}[2]{{%
\expandafter\ifx\csname l@#1\endcsname\relax
\typeout{** WARNING: IEEEtran.bst: No hyphenation pattern has been}%
\typeout{** loaded for the language `#1'. Using the pattern for}%
\typeout{** the default language instead.}%
\else
\language=\csname l@#1\endcsname
\fi
#2}}
\providecommand{\BIBdecl}{\relax}
\BIBdecl

\bibitem{KassakianSchalensee2011}
J.~Kassakian and R.~Schalensee, ``{The Future of the Electric Grid},''
  \emph{Massachusetts Institute of Technology, Technical Report}, 2011.

\bibitem{MaryB.Cain2012}
A.~C. Mary B.~Cain, Richard P.~O\rq{}Neill, ``{History of Optimal Power Flow
  and Formulations},'' \emph{Staff Report, Federal Energy Regulatory
  Commission}, 2012.

\bibitem{SadikuAlexander2011}
M.~Sadiku and C.~Alexander, \emph{{Fundamentals of Electric Circuits
  5th}}.\hskip 1em plus 0.5em minus 0.4em\relax Science Engineering \& Math,
  2011.

\bibitem{Grainger1994}
J.~J. Grainger and J.~William D.~Stevenson, \emph{{Power System
  Analysis}}.\hskip 1em plus 0.5em minus 0.4em\relax McGraw-Hill, 1994.

\bibitem{ChuaDesoerKuh1987}
L.~O. Chua, C.~A. Desoer, and E.~S. Kuh, \emph{{Linear and Nonlinear
  Circuit}}.\hskip 1em plus 0.5em minus 0.4em\relax McGraw-Hill, 1987.

\bibitem{GeneralElectric2016}
{General Electric, Siemens, V\&R POM Suite, PowerWorld, DSATools, and eTap},
  ``{Node-Breaker Modeling Representation},'' North American Electric
  Reliability Corporation (NERC), Tech. Rep., 2016.

\bibitem{ThomasKincicDaviesEtAl2016}
B.~Thomas, S.~Kincic, D.~Davies, H.~Zhang, and J.~Sanchez-Gasca, ``{A New
  Framework to Facilitate the Use of Node-Breaker Operations Model for
  Validation of Planning Dynamic Models in WECC},'' \emph{Power and Energy
  Society General Meeting (PESGM)}, 2016.

\bibitem{FischlPuntel1972}
R.~Fischl and W.~R. Puntel, ``{Computer-aided design of electric power
  transmission networks},'' \emph{IEEE PAS Conference Paper}, 1972.

\bibitem{W.R.Punteletal.1973}
W.~R.~P. et~al., ``{An automated method for long-range planning of transmission
  networks},'' \emph{Proceedings of the 10th Power Industry Computer
  Applications Conference (PICA)}, 1973.

\bibitem{HACHTELBRAYTONGUSTAVSON1972}
G.~D. HACHTEL, R.~K. BRAYTON, and F.~G. GUSTAVSON, ``{The Sparse Tableau
  Approach to Network Analysis and Design},'' \emph{IEEE Transactions on
  Circuit Theory}, vol. CT-18, no.~1, 1972.

\bibitem{NagelPederson1973}
L.~W. Nagel and D.~O. Pederson, ``{Simulation Program with Integrated Circuit
  Emphasis (SPICE)},'' \emph{16th Midwest Symp. on Circuit Theory}, 1973.

\bibitem{IBMPPDSH1973}
{IBM Program Product Document SH20-1118-0}, ``{ASTAP -- Advanced statistical
  analysis program},'' IBM Data Processing Div, Tech. Rep., 1973.

\bibitem{WEEKSJIMENEZMAHONEYEtAl1973}
W.~T. WEEKS, A.~J. JIMENEZ, G.~W. MAHONEY, D.~MEHTA, H.~QASSEMZADEH, and T.~R.
  SCOTT, ``{Algorithms for ASTAP-A Network-Analysis Program},'' \emph{IEEE
  Transactions on Circuit Theory}, vol. CT-20, no.~6, 1973.

\bibitem{DlRECTORSULLIVAN1979}
S.~W. DlRECTOR and R.~L. SULLIVAN, ``{A TABLEAU APPROACH TO POWER SYSTEM
  ANALYSIS AND DESIGN},'' \emph{Circuit Theory and Applications}, vol.~7, 1979.

\bibitem{GAMS}
\BIBentryALTinterwordspacing
GAMS, ``{General Algebraic Modeling System}.'' [Online]. Available:
  \url{https://www.gams.com/}
\BIBentrySTDinterwordspacing

\bibitem{BergenVittal2000}
A.~R. Bergen and V.~Vittal, \emph{{Power System Analysis 2nd}}.\hskip 1em plus
  0.5em minus 0.4em\relax Prentice Hall, 2000.

\bibitem{Weedy1987}
B.~M. Weedy, \emph{{Electric Power Systems 3rd}}.\hskip 1em plus 0.5em minus
  0.4em\relax John Wiley \& Sons, 1987.

\bibitem{OommenKohler1993}
M.~Oommen and J.~Kohler, ``{Effect of three-winding transformer models on the
  analysis and protection of mine power systems},'' \emph{Industry Applications
  Society Annual Meeting}, 1993.

\bibitem{Shaarbafi2014}
K.~Shaarbafi, ``{Transformer Modelling Guide},'' Alberta Electric System
  Operator (AESO), Tech. Rep., 2014.

\bibitem{LesieutreHiskens2005}
B.~C. Lesieutre and I.~A. Hiskens, ``{Convexity of the Set of Feasible
  Injections and Revenue Adequacy in FTR Markets},'' \emph{IEEE Transactions on
  Power Systems}, vol.~20, no.~4, 2005.

\bibitem{Carpentier1962}
J.~Carpentier, ``{Contribution to the Economic Dispatch Problem},'' \emph{Bull.
  Soc. Franc. Elect}, vol.~8, no.~3, 1962.

\bibitem{CongRegulskiWallEtAl2015}
Y.~Cong, P.~Regulski, P.~Wall, M.~Osborne, and V.Terzija, ``On the use of
  dynamic thermal line ratings for improving operational tripping schemes,''
  \emph{IEEE Transactions on Power Delivery}, November 2015.

\bibitem{ParkTangFerrisEtAl2016}
B.~Park, L.~Tang, M.~Ferris, and C.~L. DeMarco, ``{Examination of three
  different ACOPF formulations with generator capability curves},'' \emph{IEEE
  Transactions on Power Systems}, 2016.

\bibitem{MurliLeonePardalosEtAl1998}
A.~Murli, R.~de~Leone, P.~M. Pardalos, and G.~Toraldo, Eds., \emph{{High
  Performance Algorithms and Software in Nonlinear Optimization}}.\hskip 1em
  plus 0.5em minus 0.4em\relax Springer, 1998.

\end{thebibliography}
								
\vspace{-1cm}								
\begin{IEEEbiographynophoto}{Byungkwon Park}
received the B.S. degree in electrical engineering from the Chonbuk National University, South Korea, and the M.S degree in electrical engineering from the University of Wisconsin-Madison (UW-Madison), Madison, WI, USA, in 2011 and 2014, respectively. He is currently pursuing the Ph.D. degree in electrical engineering at UW-Madison. His research interests include optimization, electric network expansion analysis and contingency analysis of electrical energy systems.
\end{IEEEbiographynophoto}
\vspace{-1cm}
\begin{IEEEbiographynophoto}{Jayanth Netha}
received the B.Tech degree in production engineering from National Institute of Technology, Tiruchirappalli, India, in 2016. He is currently pursuing the M.S. degree in industrial engineering at UW Madison. His research interests include mathematical modeling of power systems and supply chain networks.
\end{IEEEbiographynophoto}
%\newpage
\vspace{-1cm}
\begin{IEEEbiographynophoto}{Michael C. Ferris} is the Stephen C. Kleene Professor in Computer Science, and (by courtesy) Mathematics and Industrial and Systems Engineering, and leads the Optimization Group within the Wisconsin Institutes for Discovery at the University of Wisconsin, Madison, USA. He received his PhD from the University of Cambridge, England in 1989. His research is concerned with algorithmic and interface development for large scale problems in mathematical programming, including links to the GAMS and AMPL modeling languages, and general purpose software such as PATH, NLPEC and EMP. He has worked on many applications of both optimization and complementarity, including cancer treatment planning, energy modeling, economic policy, traffic and environmental engineering, video-on-demand data delivery, structural and mechanical engineering.
\end{IEEEbiographynophoto}
\vspace{-1cm}
\begin{IEEEbiographynophoto}{Christopher L. DeMarco}
holds the Grainger Professorship in Power Engineering at the University of Wisconsin-Madison, where he been a member of the faculty of Electrical and Computer Engineering (ECE) since 1985. He has served as ECE Department Chair (2002-2005), and is UW-Madison Site Director for the Power Systems Engineering Research Center (2004-present). He was recipient of the UW-Madison Chancellor's Distinguished Teaching Award in 2000. Dr. DeMarco received his PhD degree at the University of California, Berkeley in 1985, and his B.S. degree from the Massachusetts Institute of Technology in 1980, both in Electrical Engineering and Computer Sciences. His research and teaching interests center on control, operational security, and optimization of electrical energy systems. 
\end{IEEEbiographynophoto}
\vspace{\fill}
								
%\vfill
								
								
\end{document}